\documentclass[a4paper,12pt,reqno]{amsart}

\usepackage{amsthm}
\usepackage{amsmath}
\usepackage{amssymb}
\usepackage{latexsym}
\usepackage{enumerate}
\usepackage{cite}%cite.sty
\usepackage{color}
\newtheorem{theorem}{Theorem}[section]
\newtheorem{rem}[theorem]{Remark}

\newtheorem{prop}[theorem]{Proposition}
\newtheorem{lem}[theorem]{Lemma} 
\makeatletter
\@addtoreset{equation}{section}

\makeatother
\newcommand{\N}{\mathbb{N}}
\newcommand{\Z}{\mathbb{N}_0}
\newcommand{\R}{\mathbb{R}}
    
\newcommand{\Sc}{\mathcal{S}}         
\newcommand{\Scp}{\mathcal{S}^\prime} 
\newcommand{\tp}{t^\prime}
\newcommand{\Pa}{\partial}
\newcommand{\Pt}{\partial_t}
\newcommand{\Pj}{\partial_j}
\newcommand{\Px}{\partial_x}
\newcommand{\bra}[1]{ \left<  #1 \right>  }
\newcommand{\braa}[1]{\left(  #1 \right)  }
\newcommand{\brab}[1]{\left\{ #1 \right\} }
\newcommand{\brac}[1]{\left[  #1 \right]  }
\newcommand{\abs}[1]{ \left|  #1 \right|  }
\newcommand{\nr}[1]{  \left\| #1 \right\| }
\renewcommand{\Re}{\mathrm{Re}\,}

\newcommand{\D}{{\mathsf{D}}}
\newcommand{\s}{{\mbox{\boldmath $S$}}}
\newcommand{\mh}{\mathcal{H}}
\newcommand{\oi}{\omega^{-1}}
\makeatletter %size of () ex.\bbb{3.2}\{
\newcommand{\bbb}[1]{\bBigg@{#1}} 
\makeatother
\newcommand{\Miijj}[4]{%2*2
\begin{pmatrix}
   #1 & #2 \\
   #3 & #4
\end{pmatrix}}
\DeclareMathOperator*{\supp}{supp}
\newcommand{\bib}{%

}
\begin{document}
\title[%
Scattering problems for NLKG]{%
Weighted Sobolev regularity and failure of profile convergence 
associated with the scattering operator 
for the three-dimensional cubic Klein-Gordon equation}

\author[H. Sasaki]{%
Hironobu Sasaki$^*$\\
Department of Mathematics and Informatics, 
Chiba University, 263--8522, Japan.}

\thanks{$^*$Supported by 
JSPS KAKENHI Grant Number 22K03367.}

\thanks{$^*$email:
\texttt{sasaki@math.s.chiba-u.ac.jp}} 
\subjclass[2020]{Primary 35P25; Secondary 35B40, 35L71.}
%MSC2020
%35P25   	Scattering theory for PDEs

\keywords{%
Nonlinear Klein-Gordon equation; 
Scattering problem}

\begin{abstract}%\subsubsection{abstract}
The scattering problem 
for the three-dimensional nonlinear Klein-Gordon equation 
\begin{align*}
u_{tt}-\Delta u+u=\pm u^3
\end{align*} 
is studied. 
It was shown that 
the scattering operator ${\mbox{\boldmath $S$}}$ 
for the equation is well-defined 
on a neighborhood of 0 in the Sobolev space 
$\mathcal{H}^{1/2}:=H^{1/2}(\mathbb{R}^3) \oplus H^{-1/2}(\mathbb{R}^3)$, 
and that 
if $\Phi_-$ is sufficiently small in $\mathcal{H}^{1/2}$ 
and belongs to $\mathcal{S} \oplus \mathcal{S}$,   
then 
${\mbox{\boldmath $S$}}\Phi_-$ is also in $\mathcal{S} \oplus \mathcal{S}$, 
where 
$\mathcal{S}=\mathcal{S}(\R^3)$ is the Schwartz space. 
In the present paper, however, we prove that
for any non-zero real-valued $\Phi_-\in\mathcal{S}\oplus\mathcal{S}$
whose $\mathcal{H}^{1/2}$-norm is sufficiently small,
the corresponding solution $u$ satisfies
\[
U(-t) \,{}^t(u(t),\partial_t u(t))
\not\to {\mbox{\boldmath $S$}}\Phi_- 
\quad(t\to+\infty) 
\]
in the weighted Sobolev space 
$H^{1,3}(\mathbb{R}^3)\oplus H^{0,3}(\mathbb{R}^3)$ 
although $U(-t)\,{}^t(u(t),\partial_t u(t))$
converges to ${\mbox{\boldmath $S$}}\Phi_-$
in $H^{1,\sigma}(\mathbb{R}^3)\oplus H^{0,\sigma}(\mathbb{R}^3)$
for any $\sigma<3$.
Here, $U(t)$ denotes the free Klein-Gordon evolution group.
\end{abstract}

\maketitle

\section{Introduction}
We consider 
the scattering problem 
for the nonlinear Klein-Gordon equation 
\begin{align}\label{NLKG}
\Pt^2 u-\Delta u+u=\lambda u^3
\end{align} 
in three space dimensions. 
Here, 
$u=u(t,x)$ is a real-valued unknown function 
defined on $\R \times \R^3$, 
$\Pt=\Pa/\Pa t$, 
$\Delta=\Pa_1^2+\Pa_2^2+\Pa_3^2$ is the Laplacian in $x$, 
$\Pj=\Pa/\Pa x_j$ ($j=1,2,3$) 
and $\lambda$ is a non-zero real constant.
By rescaling $u$, 
we may assume without loss of generality
that $\lambda\in\{-1,1\}$.
There is a large literature on the scattering problem 
for nonlinear Klein-Gordon equations 
(see, e.g., \cite{%
Brenner1982,
Brenner1985,
GinibreVelo1985,
HayashiNaumkin2008,
HayashiNaumkin2009,
IbrahimMasmoudiNakanishi2009,
IkedaInuiOkamoto2020,
MasakiSegata2018,
MasakiSegata2018-2,
NakamuraOzawa2001,
NakamuraOzawa2002,
Sasaki2020,
Sasaki2024,
Wang1998} 
and references therein). 
In particular,
Wang \cite{Wang1998} proved that
the scattering operator $\s$ for Equation (\ref{NLKG})
is well-defined on a neighborhood of $0$
in an appropriate function space.

To state this existence result and the aim of this paper precisely,
we introduce some notation.
For $p\in [1,\infty]$, 
$L^p$ and $\nr{\cdot}_p$ denote $L^p(\R^3)$ and its norm, 
respectively.   
Set 
$\bra{x}=\sqrt{1+|x|^2}$ ($x\in \R^3$) 
and 
$\omega=\sqrt{1-\Delta}$. 
For $s,\sigma\in \R$, 
we define 
the weighted Sobolev space 
$H^{s,\sigma}=\omega^{-s}\bra{x}^{-\sigma}L^2(\R^3)$
and 
$\mh^{s,\sigma}=H^{s,\sigma}\oplus H^{s-1,\sigma}$. 
For $\Phi={}^t (f,g) \in \mh^{s,\sigma}$, 
the $\mh^{s,\sigma}$-norm is given by 
\[
\nr{\Phi}_{\mh^{s,\sigma}}=
\sqrt{
\nr{\bra{x}^\sigma \omega^s f}_2^2 
+ 
\nr{\bra{x}^\sigma \omega^{s-1}g}_2^2
}.
\]
We write $H^s=H^{s,0}$ and $\mh^s=\mh^{s,0}$.
We denote the Schwartz space 
(resp. the set of all tempered distributions) on $\R^3$   
by $\Sc$ (resp. $\Scp$). 
Put 
\begin{align*}
U(t)=
\Miijj{\cos(t\omega)}{\oi\sin(t\omega)}{-\omega\sin(t\omega)}{\cos(t\omega)},
\quad
t\in \R.
\end{align*} 
For any $s\in\R$, $(U(t))_{t\in\R}$ is
a strongly continuous one-parameter unitary group on $\mh^s$.
For $\Phi={}^t(f,g)\in \Scp \oplus \Scp$, 
we define 
\begin{align*}
w[\Phi](t)=\cos(t\omega)f+\oi\sin(t\omega)g,
\quad
t\in \R, 
\end{align*} 
which solves the free Klein-Gordon equation
$\Pt^2u-\Delta u+u=0$
with initial conditions $u(0)=f$ and $\Pt u(0)=g$
in the sense of distributions.
The operator $\D$ is given by $\D={}^t(1,\Pt)$.
Then we have 
\begin{align*}
U(t)\Phi=
\D w[\Phi](t)
=
\binom{
w[\Phi](t)
}{
\partial_t w[\Phi](t)
},
\quad
t\in \R.
\end{align*} 
Therefore, 
$U(t)$ is called the free Klein-Gordon evolution group. 
For a Banach space $X$ equipped with the norm $\nr{\cdot}$ 
and $\eta>0$,  
we denote the closed ball 
$\{ f\in X \,;\, \nr{f}\le \eta \}$
by $B_\eta X$. 
Furthermore, for $p\in[1,\infty]$, we write
$L^pX=L^p(\R;X)$ and define
\[
C_bX=C(\R;X)\cap L^\infty X,
\qquad
C_b^1X=
\brab{u\in C^1(\R;X)\,;\,u, \Pt u \in C_bX}.
\]
Define 
$Z=C_b H^{1/2} \cap C_b^1 H^{-1/2} \cap L^4L^4$.
With this notation, 
the existence result mentioned above
can be stated as follows:

\begin{prop}\label{prop:S}%\paragraph{prop:S}
There exists an $\eta_1>0$ such that,
for any $\eta\in(0,\eta_1]$ and any
$\Phi_-={}^t(f_-,g_-)\in B_\eta\mh^{1/2}$,
there exist a unique global solution 
$u \in Z$  
satisfying
\begin{align}\label{Ieq}
u(t)
=
w[\Phi_-](t)
+\lambda 
\int_{-\infty}^t
\oi \sin((t-\tp)\omega) u(\tp)^3 d\tp 
\ \text{in $H^{1/2}$,} 
\qquad t\in\R,
\end{align}
and a unique state $\Phi_+\in\mh^{1/2}$ satisfying
\begin{align}\label{conv}
\lim_{t\to\pm\infty}U(-t)\D u(t)=\Phi_\pm
\end{align}
in $\mh^{1/2}$.
Moreover, $\Phi_+\in B_{2\eta}\mh^{1/2}$.
Hence the scattering operator 
$\s:B_{\eta_1}\mh^{1/2}\ni \Phi_- \mapsto \Phi_+ \in B_{2\eta_1}\mh^{1/2}$ 
is well-defined. 
By decreasing $\eta_1$ if necessary,
we may assume that $\s$ is injective.
\end{prop}

We call $U(-t)\D u(t)$ the profile of the solution $u$. 
In this paper, 
we consider the role of the weighted Sobolev spaces 
$\mh^{s,\sigma}$ in the study of $\s$.
For some nonlinear Klein-Gordon equations and related equations
with ``$L^2$-subcritical'' nonlinearities,
the spaces $\mh^{s,1}$ play an important role
in proving the existence of scattering operators.
For example,
for the three-dimensional Klein-Gordon equation with nonlinearity
$\pm |u|^{p-1}u$, where $9/5<p<7/3$, 
the scattering operator is well-defined
on a neighborhood of $0$ in $\mh^{5/3,1}$, 
and 
the convergence (\ref{conv}) holds in the same space
(see Hayashi and Naumkin 
\cite{HayashiNaumkin2008} and also \cite{HayashiNaumkin2009}). 
In \cite{HayashiNaumkin2008},
the linear operators
$J_\pm=e^{\pm it\omega}xe^{\mp it\omega}$
were used to obtain time-decay estimates,
while the closely related Lorentz boost operators
\begin{align}\label{def:P}
P_j=t\Pj+x_j\Pt \quad (j=1,2,3)
\end{align}
were employed to estimate the nonlinear term. 
Recently, 
Yang \cite{Yang2026} considered 
the two-dimensional semi-relativistic Hartree equation 
$iu_t-\omega u = \lambda (|x|^{-\gamma} \ast |u|^2)u$  
with $1<\gamma<2$,
and proved scattering in $H^s(\R^2)\cap H^{10,2}(\R^2)$
for sufficiently large $s$ and sufficiently small
initial data in the same space.

Returning to (\ref{NLKG}),
we see from Proposition \ref{prop:S} that
no positive spatial weight is needed
to establish the existence of $\s$. 
However, 
it is natural to ask whether $\s$ 
preserves regularity and decay properties,
and whether the convergence (\ref{conv})
also holds in the corresponding weighted Sobolev spaces.
The former question was addressed in \cite{Sasaki2020}.
However, these preservation results do not determine
whether the convergence (\ref{conv}) holds in the
corresponding weighted Sobolev spaces.
The argument in \cite{Sasaki2020} shows that 
$\s$ preserves $\mh^{N+1/2,N}$ for every $N\in \N$,
provided that the input data $\Phi_-$ 
is sufficiently small in $\mh^{1/2}$. 
In particular, 
$\s$ preserves the Schwartz space
$\Sc\oplus\Sc$. 
See also Simon and Taflin \cite{SimonTaflin1993}
for related results under stronger smallness assumptions.
For nonlinear Schr\"odinger equations,
Galilean vector fields provide a standard framework
for studying weighted regularity and decay
(see, e.g., \cite{HayashiNakamitsuTsutsumi1988}).
One might therefore expect the convergence (\ref{conv})
to hold also in $\mh^{N+1/2,N}$. 
The main result of the present paper, however,
shows that this expectation is false in general:

\begin{theorem}\label{thm:main}%\paragraph{thm:main}
Let $\eta_1$ and $\s$ be the positive number and the mapping
appearing in Proposition \ref{prop:S}, respectively.
Then there exists an $\eta_2\in(0,\eta_1]$ such that,
for any $\eta\in(0,\eta_2]$ and any non-zero real-valued
$\Phi_-\in B_\eta\mh^{1/2}\cap(\Sc\oplus\Sc)$,
we have
$\s\Phi_-\in B_{2\eta}\mh^{1/2}\cap(\Sc\oplus\Sc)$,
and the corresponding solution $u$ to (\ref{Ieq}) satisfies
$U(-t)\D u(t)\in\Sc\oplus\Sc$ for any $t\in\R$.
Moreover,
\begin{gather}
\label{conv:1}
\lim_{t\to+\infty}
\nr{U(-t)\D u(t)-\s\Phi_-}_{\mh^{1,\sigma}}
=0
\quad (0\le\sigma<3),
\\
\label{core}
\begin{aligned}
0
  <\liminf_{t\to+\infty}
  \nr{U(-t)\D u(t)-\s\Phi_-}_{\mh^{1,3}}
  \le\limsup_{t\to+\infty}
  \nr{U(-t)\D u(t)-\s\Phi_-}_{\mh^{1,3}}
<\infty,
\end{aligned}
\\
\label{conv:2}
\lim_{t\to+\infty}
\nr{U(-t)\D u(t)-\s\Phi_-}_{\mh^{1,\sigma}}
=\infty
\quad (\sigma>3).
\end{gather}
\end{theorem}

\begin{rem}\label{rem:1:thm:main}%\paragraph{rem:1:thm:main}
Theorem \ref{thm:main} identifies $\sigma=3$
as the sharp threshold for convergence of
the profile $U(-t)\D u(t)$ to $\s\Phi_-$
in $\mh^{1,\sigma}$ as $t\to+\infty$.
At the threshold, there exist constants $c,C>0$ such that
\[
c\le
\nr{U(-t)\D u(t)-\s\Phi_-}_{\mh^{1,3}}
\le C
\]
for all sufficiently large $t$.
For any $\sigma>3$, the norm of the difference
in $\mh^{1,\sigma}$ tends to infinity as $t\to+\infty$.
\end{rem}

\begin{rem}\label{rem:2:thm:main}%\paragraph{rem:2:thm:main}
The same conclusions hold for the half-wave profiles.
Indeed, set
\[
u_\pm(t)=\frac12\braa{u(t)\mp i\oi\Pt u(t)},
\qquad
\phi_\pm=\frac12\braa{f_+\mp i\oi g_+},
\]
where $\s\Phi_-={}^t(f_+,g_+)$.
Since $u$, $f_+$ and $g_+$ are real-valued, we have
\[
\nr{e^{\mp it\omega}u_\pm(t)-\phi_\pm}_{H^{1,\sigma}}
=
\frac{1}{2}
\nr{U(-t)\D u(t)-\s\Phi_-}_{\mh^{1,\sigma}}
\quad (\sigma\in \R).
\]
Thus the convergence threshold is unchanged
under this decomposition.
\end{rem}

%\subsection{Strategy}
We now describe the strategy for proving
the first inequality in (\ref{core}),
which implies the failure of profile convergence. 
For non-negative numbers $a$ and $b$, 
we write $a\lesssim b$ or $b\gtrsim a$ if $a\le Cb$
for some constant $C>0$ 
independent of the time and space variables.
The constant may vary from line to line and may depend
on fixed parameters and the fixed initial data.
Let $\Phi_-$ and $u$ be
the data and solution appearing in Theorem \ref{thm:main}. 
Set  
\begin{align}\label{def:h}
h(t)=U(-t)\D u(t)-\s\Phi_-,
\quad
t\in \R.
\end{align} 
We choose some $j\in\{1,2,3\}$,
apply the Lorentz boost operator $P:=P_j$ three times to $u$,
and use the operators
\begin{align}\label{def:AB}
A(t)
:=\Miijj{t\Pj}{x_j}{\Pj- x_j \omega^2}{t\Pj}
\quad\text{and}\quad
B:=A(0)
=\Miijj{0}{x_j}{-\omega x_j \omega}{0}
=\Miijj{0}{x_j}{\Pj- x_j \omega^2}{0},
\end{align} 
which are related to $P$ 
through the identities 
\begin{align}\label{AUUB}
A(t)U(t)=U(t)B
\quad\text{and}\quad
\D Pu(t)
=
A(t)\D u(t)+\binom{0}{\lambda x_j u^3}.
\end{align}
Furthermore, 
note that $B^3$ is bounded from $\mh^{1,3}$ to $\mh^{-2}$. 
We see from the unitarity of $U(t)$ on $\mh^{-2}$ that 
\begin{align*}
\nr{h(t)}_{\mh^{1,3}}
\gtrsim 
\nr{B^3 h(t)}_{\mh^{-2}}
=
\nr{U(t)B^3 h(t)}_{\mh^{-2}}
\ge 
\nr{(U(t)B^3 h(t))_1}_{H^{-2}}.
\end{align*} 
By (\ref{AUUB}), 
we obtain 
\begin{align}
(U(t)B^3 h(t))_1
&=
P^3(u-u_+)(t)-3tx_jF - 2tx_j^2\Pj (F-F_+) - x_j^2 PF - 2tx_j^2\Pj F_+ 
\nonumber 
\\
&=:
d(t) - 2tx_j^2\Pj F_+,
\label{def:d}
\end{align} 
where
$u_+:=w[\s \Phi_-]$, 
$F=\lambda u^3$ and $F_+=\lambda u_+^3$. 
The key to the proof is to choose a family of test functions
$\psi_t$ on $\R^3$ such that 
\begin{align}
\limsup_{t\to\infty}\nr{\psi_t}_{H^2}
&<\infty,
\label{property:psi:1}
\\
\bra{d(t),\psi_t}_{L^2}
&=o(1) \quad\text{as $t\to\infty$}
\label{property:psi:2}
\end{align} 
and 
\begin{align}
\bra{tx_j^2\Pj F_+,\psi_t}_{L^2}
=C_\psi+O(t^{-1})
\quad\text{as $t\to\infty$}.
\label{property:psi:3}
\end{align} 
Here, 
$C_\psi$ is a positive constant. 
These properties imply the first inequality in (\ref{core}).
The stationary phase expansion of $u_+$ allows us
to construct test functions $\psi_t$
that detect the nonvanishing contribution of $tx_j^2\Pj F_+$
and whose $H^2$-norms are uniformly bounded in $t$.
\\

The rest of this paper is organized as follows.
Section 2 collects some preliminary results.
In Section 3, we establish the regularity and decay properties of the solution.
Section 4 is devoted to the construction of the test functions.
Finally, we prove Theorem \ref{thm:main} in Section 5.

\section{Preliminaries}
In this section, 
we introduce notation and propositions 
used in the rest of this paper. 
For non-negative numbers $a$ and $b$, 
we write $a\sim b$ if $a\lesssim b$ and $b\lesssim a$.
For real numbers $a$ and $b$, 
set $a\vee b=\max\{a,b\}$
and $a\wedge b=\min\{a,b\}$.
We denote the set of all non-negative integers by $\Z$. 
For $f\in \Sc$,  
the Fourier transform of $f$ is given by 
\[
\widehat f(\xi)
=(2\pi)^{-3/2}\int_{\R^3}e^{-ix\cdot\xi}f(x)\,dx  
\quad (\xi \in \R^3).
\]
For $k\in\Z$ and $1\le p\le\infty$, 
we define the classical Sobolev space $W^k_p$
and its norm $\nr{\cdot}_{W^k_p}$ by 
\[
W^k_p=\{f\in L^p;\ \partial_x^\alpha f\in L^p
\text{ for }|\alpha|\le k\},
\qquad
\|f\|_{W^k_p}=\sum_{|\alpha|\le k}\|\partial_x^\alpha f\|_p.
\]
Here, 
the derivatives are understood in the sense of distributions. 
Put $Y=L^4L^4$. 
We equip $Z$ with the norm
\[
\|v\|_Z=
\|v\|_{L^\infty H^{1/2}}
\vee\|\partial_t v\|_{L^\infty H^{-1/2}}
\vee\|v\|_Y.
\]
For a suitable function $G$, define
\[
\mathcal I_s[G](t)
=\int_{-\infty}^t
\omega^{-1}\sin((t-\tau)\omega)G(\tau)\,d\tau,
\qquad
\mathcal I_c[G](t)
=\int_{-\infty}^t
\cos((t-\tau)\omega)G(\tau)\,d\tau.
\]
We also define the Duhamel operators with integration
from $t$ to $+\infty$ by 
\[
\mathcal I_s^+[G](t)
=-\int_t^\infty
\omega^{-1}\sin((t-\tau)\omega)G(\tau)\,d\tau,
\qquad
\mathcal I_c^+[G](t)
=-\int_t^\infty
\cos((t-\tau)\omega)G(\tau)\,d\tau.
\]
Thus $\partial_t\mathcal I_s[G]=\mathcal I_c[G]$
and $\partial_t\mathcal I_s^+[G]=\mathcal I_c^+[G]$.

\begin{prop}\label{prop:free-decay}%\paragraph{prop:free-decay}
There exists a constant $C>0$ such that
\[
\|w[\Phi](t)\|_\infty
\le C|t|^{-3/2}
\bigl(\|f\|_{W^3_1}+\|g\|_{W^2_1}\bigr)
\qquad (t\ne0)
\]
for $\Phi={}^t(f,g)\in W^3_1\oplus W^2_1$.
Moreover,
\[
\|w[\Phi](t)\|_\infty
\le C\bigl(\|f\|_{H^2}+\|g\|_{H^1}\bigr)
\qquad (t\in\R)
\]
for $\Phi\in\mh^2$.
\end{prop}

\begin{proof}
Let $B^s_{p,q}$ denote the inhomogeneous Besov space
on $\R^3$, and fix $0<a<1/2$.
By \cite[(2.12)]{NakamuraOzawa2002}, with
$n=3$, $\lambda_1=5/2$ and $\sigma_1=3$,
and the lifting property of Besov spaces, we have
\[
\|e^{\pm it\omega}h\|_{B^a_{\infty,2}}
\lesssim |t|^{-3/2}\|h\|_{B^{5/2+a}_{1,2}}
\qquad (t\ne0).
\]
Since
\[
B^a_{\infty,2}\hookrightarrow L^\infty,
\qquad
W^3_1\hookrightarrow B^{5/2+a}_{1,2},
\qquad
W^2_1\hookrightarrow B^{3/2+a}_{1,2},
\]
the representation
\[
w[\Phi](t)=\cos(t\omega)f+\omega^{-1}\sin(t\omega)g
\]
gives
\begin{align*}
\|w[\Phi](t)\|_\infty
&\lesssim |t|^{-3/2}
\left(
\|f\|_{B^{5/2+a}_{1,2}}
+\|\omega^{-1}g\|_{B^{5/2+a}_{1,2}}
\right)
\\
&\lesssim |t|^{-3/2}
\left(\|f\|_{W^3_1}+\|g\|_{W^2_1}\right).
\end{align*}
The second estimate follows from the embedding
$H^2\hookrightarrow L^\infty$ and
\[
\|w[\Phi](t)\|_{H^2}
\le \|f\|_{H^2}+\|g\|_{H^1}.
\]
\end{proof}

\begin{prop}\label{prop:St}%\paragraph{prop:St}
There exists a constant $C>0$ such that
\[
\|w[\Phi]\|_Z\le C\|\Phi\|_{\mh^{1/2}}
\]
for any $\Phi\in\mh^{1/2}$, and
\[
\|\mathcal I_s[G]\|_Z+
\|\mathcal I_s^+[G]\|_Z
\le C\|G\|_{L^{4/3}L^{4/3}}
\]
for any $G\in L^{4/3}L^{4/3}$.
In particular,
\[
\left\|
\int_t^\infty U(t-\tau)\binom{0}{G(\tau)}\,d\tau
\right\|_{\mh^{1/2}}
\le C\|G\|_{L^{4/3}((t,\infty);L^{4/3})}.
\]
The Duhamel operators are understood as the continuous
extensions of the operators initially defined on
$C_c^\infty(\R^{1+3})$.
\end{prop}

\begin{proof}
The assertions follow from
\cite[Lemma 2.1]{NakamuraOzawa2001},
the embeddings
$B^0_{4,2}\hookrightarrow L^4$ and
$L^{4/3}\hookrightarrow B^0_{4/3,2}$,
and time reversal.
\end{proof}

\begin{prop}\label{prop:AUUB}%\paragraph{prop:AUUB}
For $j=1,2,3$, define
\[
A_j(t)=
\Miijj{t\partial_j}{x_j}
{\partial_j-x_j\omega^2}{t\partial_j},
\qquad
B_j=A_j(0).
\]
Then
\[
A_j(t)U(t)=U(t)B_j.
\]
If $v$ satisfies
\[
(\partial_t^2-\Delta+1)v=G,
\]
then
\[
\D(P_jv)=A_j(t)\D v+\binom{0}{x_jG}.
\]
Consequently, for any $n\ge1$,
\[
\D(P_j^nv)
=A_j(t)^n\D v+
\sum_{k=0}^{n-1}
A_j(t)^{n-1-k}\binom{0}{x_jP_j^kG}.
\]
These identities hold for smooth functions.
They extend to distributions by approximation,
provided that all terms are well-defined.
\end{prop}

\begin{proof}
The identity
$\partial_j-x_j\omega^2=-\omega x_j\omega$
follows from $[\omega,x_j]=-\omega^{-1}\partial_j$.
Direct differentiation gives
\[
\partial_t(P_jv)
=t\partial_j\partial_t v+\partial_jv
-x_j\omega^2v+x_jG,
\]
which proves the second identity.
For a free solution $v=w[\Phi]$, $P_jv$ is again a free
solution and has Cauchy data $B_j\Phi$ at $t=0$.
Hence
\[
A_j(t)U(t)\Phi=\D(P_jw[\Phi])(t)=U(t)B_j\Phi.
\]
The last identity follows by induction, since $P_j$
commutes with $\partial_t^2-\Delta+1$.
\end{proof}

\begin{prop}\label{prop:B-bdd}%\paragraph{prop:B-bdd}
There exists a constant $C\ge1$ such that
\[
C^{-1}\|\Phi\|_{\mh^{1,3}}
\le
\|\Phi\|_{\mh^1}
+\sum_{j=1}^3\|B_j^3\Phi\|_{\mh^{-2}}
\le
C\|\Phi\|_{\mh^{1,3}}
\]
for any $\Phi\in\Sc\oplus\Sc$.
\end{prop}

\begin{proof}
Set
\[
T=\begin{pmatrix}\omega&0\\0&1\end{pmatrix},
\qquad
J=\begin{pmatrix}0&1\\-1&0\end{pmatrix},
\qquad q=T\Phi.
\]
Since $TB_jT^{-1}=J\omega x_j$, we have
\[
\|\Phi\|_{\mh^1}=\|q\|_2,\qquad
\|\Phi\|_{\mh^{1,3}}=\|\langle x\rangle^3q\|_2,
\qquad
\|B_j^3\Phi\|_{\mh^{-2}}
=\|\omega^{-3}(\omega x_j)^3q\|_2.
\]
By \cite[Proposition 2.3]{Sasaki2020} and the
commutation relation between $x_j$ and Fourier multipliers,
we can write
\[
\omega^{-3}(\omega x_j)^3
=x_j^3+\sum_{k=0}^2M_{j,k}x_j^k,
\]
where each $M_{j,k}$ is a bounded Fourier multiplier on $L^2$.
Consequently,
\[
\sum_{j=1}^3
\|\omega^{-3}(\omega x_j)^3q-x_j^3q\|_2
\lesssim \|\langle x\rangle^2q\|_2.
\]
Using
\[
\|\langle x\rangle^3q\|_2
\sim \|q\|_2+\sum_{j=1}^3\|x_j^3q\|_2
\]
and, for any $\varepsilon>0$,
\[
\|\langle x\rangle^2q\|_2
\le \varepsilon\|\langle x\rangle^3q\|_2
+C_\varepsilon\|q\|_2,
\]
we obtain both inequalities by choosing $\varepsilon$
sufficiently small.
\end{proof}

\begin{prop}\label{prop:PI}%\paragraph{prop:PI}
Let $j\in\{1,2,3\}$ and $n\in\Z$.
Suppose that
\[
P_j^kG\in L^{4/3}(\R;L^{4/3})
\qquad (0\le k\le n),
\]
where the derivatives are understood in the sense
of distributions.
Then
\[
P_j^n\mathcal I_s[G]=\mathcal I_s[P_j^nG],
\qquad
P_j^n\mathcal I_s^+[G]=\mathcal I_s^+[P_j^nG]
\]
in the sense of distributions.
Spatial derivatives also commute with these operators,
provided that the corresponding derivatives of $G$
belong to $L^{4/3}L^{4/3}$.
\end{prop}

\begin{proof}
Write $p=4/3$ and $P=P_j$.
We first prove the assertion for $n=1$.
For $G\in C_c^\infty(\R^{1+3})$, the identity
\[
P\mathcal I_s[G]=\mathcal I_s[PG]
\]
follows from \cite[Proposition 3.7]{Sasaki2020}.
Now suppose that $G,PG\in L^p(\R^{1+3})$.
By the same cutoff and mollification argument as in
the proof of \cite[Proposition A.4]{Sasaki2020},
with $L^4L^4$ replaced by $L^p(\R^{1+3})$,
there exists a sequence
$(G_\ell)_{\ell\ge1}\subset C_c^\infty(\R^{1+3})$
such that
\[
G_\ell\to G,\qquad PG_\ell\to PG
\quad\text{in }L^p(\R^{1+3}).
\]

By Proposition \ref{prop:St}, we have
\[
\mathcal I_s[G_\ell]\to\mathcal I_s[G],
\qquad
\mathcal I_s[PG_\ell]\to\mathcal I_s[PG]
\quad\text{in }Z.
\]
Passing to the limit in
$P\mathcal I_s[G_\ell]=\mathcal I_s[PG_\ell]$
in the sense of distributions yields
\[
P\mathcal I_s[G]=\mathcal I_s[PG].
\]
Iteration proves the assertion for arbitrary $n$.

To treat $\mathcal I_s^+$, define
$(\mathcal R G)(t,x)=G(-t,x)$.
Then
\[
\mathcal I_s^+=\mathcal R\mathcal I_s\mathcal R,
\qquad
P\mathcal R=-\mathcal R P,
\]
so the corresponding identity follows by time reversal.
Finally, the assertions for spatial derivatives follow
from the translation invariance of the Duhamel operators.
\end{proof}

\section{Regularity and decay properties}
In the proof of the main theorem,
we need certain regularity and decay properties of the solution $u$.
For this purpose, we establish the following propositions.

\begin{prop}\label{prop:r-d}%\paragraph{prop:r-d}
Let $\eta_1$ and $\s$ be the positive number and the mapping
appearing in Proposition \ref{prop:S}, respectively. 
There exists some $\eta_2 \in (0,\eta_1]$ such that 
if $\eta \in (0,\eta_2]$ 
and $\Phi_- \in B_\eta \mh^{1/2} \cap (\Sc \oplus \Sc)$, 
then the corresponding solution $u$ to (\ref{Ieq}) satisfies
\begin{align}\label{PPu}
\Px^\alpha P_j^n u\in Z
\quad
(\alpha\in\Z^3,\ n\in\Z,\ j=1,2,3)
\end{align} 
and 
\begin{align}\label{SinS}
\s\Phi_-\in B_{2\eta}\mh^{1/2}\cap(\Sc\oplus\Sc).
\end{align} 
\end{prop}

\begin{proof}
By decreasing $\eta_2$ if necessary,
the property (\ref{SinS}) follows directly from
\cite[Theorem 1.3(4)]{Sasaki2020}.
The property (\ref{PPu}) follows by a similar argument,
which we outline for completeness.
Fix $j\in\{1,2,3\}$, $\eta>0$, and
$\Phi_-\in B_\eta\mh^{1/2}\cap(\Sc\oplus\Sc)$.
We may assume that $\Phi_-\neq 0$. 
Write $P=P_j$, $A(t)=A_j(t)$, and $B=B_j$.
Define 
the zero multi-index $o=(0,0,0)$, 
$\bra{\alpha}=1+|\alpha|$ ($\alpha \in \Z^3$),  
and the function space 
\[
Z_\infty=\brab{
v\in Z \,;\, \Px^\alpha P^n v \in Z \ (\alpha\in \Z^3, \ n\in \Z)
}.
\] 
Fix $\alpha\in \Z^3$, $n\in \Z$ and $v\in Z_\infty$. 
By Propositions \ref{prop:AUUB} and  \ref{prop:St}, 
we have $w[\Phi_-] \in Z_\infty$ 
and 
\begin{align*}
\nr{ \Px^\alpha P^n w[\Phi_-] }_Z 
\lesssim \nr{ \Px^\alpha B^n \Phi_-}_{\mh^{1/2}}. 
\end{align*} 
Here and below, 
the implicit constants are independent of $\alpha$ and $n$.
We see from Proposition \ref{prop:PI} and the Leibniz rule that 
\begin{align*}
\Px^\alpha P^n \mathcal{I}_s[v^3]
=
\sum_{\beta+\gamma+\delta=\alpha}
\sum_{m_1+m_2+m_3        =n}
\frac{\alpha!}{\beta!\gamma!\delta!}
\frac{n!}{m_1!m_2!m_3!}
\mathcal{I}_s \brac{
(\Px^\beta  P^{m_1}v) 
(\Px^\gamma P^{m_2}v) 
(\Px^\delta P^{m_3}v) 
}
\end{align*} 
and 
\begin{align*}
\Pt \Px^\alpha P^n \mathcal{I}_s[v^3]
=
\sum_{\beta+\gamma+\delta=\alpha}
\sum_{m_1+m_2+m_3        =n}
\frac{\alpha!}{\beta!\gamma!\delta!}
\frac{n!}{m_1!m_2!m_3!}
\mathcal{I}_c \brac{
(\Px^\beta  P^{m_1}v) 
(\Px^\gamma P^{m_2}v) 
(\Px^\delta P^{m_3}v) 
}.
\end{align*}
Furthermore, 
Proposition \ref{prop:St} and H\"older's inequality imply that 
\begin{align*}
&
\nr{\Px^\alpha P^n \mathcal{I}_s[v^3]}_Z
\\
&\lesssim
\sum_{\beta+\gamma+\delta=\alpha}
\sum_{m_1+m_2+m_3        =n}
\frac{\alpha!}{\beta!\gamma!\delta!}
\frac{n!}{m_1!m_2!m_3!}
\nr{\Px^\beta  P^{m_1}v}_Y 
\nr{\Px^\gamma P^{m_2}v}_Y 
\nr{\Px^\delta P^{m_3}v}_Y.
\end{align*} 
Hence the function 
$\widetilde{v}:=w[\Phi_-] +\lambda \mathcal{I}_s[v^3]$ 
belongs to $Z_\infty$ and satisfies  
\begin{align}\label{Z-infty}
&
\nr{\Px^\alpha P^n \widetilde{v}}_Z
\lesssim
\nr{\Px^\alpha B^n \Phi_-}_{\mh^{1/2}}
\nonumber
\\ \quad 
&+
\sum_{\beta+\gamma+\delta=\alpha}
\sum_{m_1+m_2+m_3        =n}
\frac{\alpha!}{\beta!\gamma!\delta!}
\frac{n!}{m_1!m_2!m_3!}
\nr{\Px^\beta  P^{m_1}v}_Y 
\nr{\Px^\gamma P^{m_2}v}_Y 
\nr{\Px^\delta P^{m_3}v}_Y.
\end{align}

We define $r(o,0)=1$ and
\begin{align*}
r(\beta,m)
=
\left\{
\max_{\substack{
o\le\gamma\le\beta,\ 0\le\ell\le m\\
|\gamma|+\ell>0
}}
\left(
\frac{
\bra{\gamma}^4(1+\ell)^2
\nr{\Px^\gamma B^\ell\Phi_-}_{\mh^{1/2}}
}{
\gamma!\ell!\nr{\Phi_-}_{\mh^{1/2}}
}
\right)^{1/(|\gamma|+\ell)}
\right\}^{|\beta|+m}
\end{align*}
for $\beta\in\Z^3$ and $m\in\Z$ with $|\beta|+m>0$.
We have 
\[
r(\beta,m_1) r(\gamma,m_2) r(\delta,m_3) \le r(\alpha,n)
\]
for any 
$\beta, \gamma, \delta \in \Z^3$ 
and 
$m_1,m_2,m_3 \in \Z$ 
with 
$\beta+\gamma+\delta=\alpha$
and 
$m_1+m_2+m_3        =n$.  
Moreover, 
the number 
\begin{align*}
\sup_{\alpha\in\Z^3,\ n\in\Z}
\sum_{\beta+\gamma+\delta=\alpha}
\sum_{m_1+m_2+m_3        =n}
\braa{ \frac{ \bra{\alpha} }{ \bra{\beta} \bra{\gamma} \bra{\delta} } }^4
\braa{ \frac{(1+n)}{(1+m_1)(1+m_2)(1+m_3)} }^2 
\end{align*} 
is finite 
(see \cite[Proposition 2.4(1)]{Sasaki2016}
with $(d,p)=(3,4)$ and $(1,2)$). 
Combining these facts with (\ref{Z-infty}), 
we obtain  
\begin{align*}
&
\frac{\bra{\alpha}^4 (1+n)^2 \nr{\Px^\alpha P^n \widetilde{v}}_Z}{
\alpha! n! r(\alpha,n)}
\\
&\lesssim
\frac{\bra{\alpha}^4 (1+n)^2 \nr{\Px^\alpha B^n \Phi_-}_{\mh^{1/2}}}{
\alpha! n! r(\alpha,n)}
+ 
\braa{
\max_{o\le \beta\le \alpha, \ 0\le m\le n}
\frac{\bra{\beta}^4 (1+m)^2 \nr{\Px^\beta P^{m} v}_Z}{
\beta! m! r(\beta,m)}
}^3 .
\end{align*} 
Therefore, 
if we define the function space 
\begin{align*}
Z[\Phi_-]
=
\brab{
v\in Z_\infty \,;\,
\nr{v}_{Z[\Phi_-]} :=
\sup_{\alpha\in \Z^3, \ n\in \Z} 
\frac{\bra{\alpha}^4 (1+n)^2 \nr{\Px^\alpha P^n v}_Z}{
\alpha! n! r(\alpha,n)}
<\infty
},
\end{align*} 
then for any $v\in Z[\Phi_-]$ we have 
$\widetilde{v}\in  Z[\Phi_-]$ and  
\begin{align*}
\nr{\widetilde{v}}_{Z[\Phi_-]}
\lesssim
\nr{\Phi_-}_{\mh^{1/2}}+\nr{v}_{Z[\Phi_-]}^3.
\end{align*} 
Similarly, for any $v,w\in Z[\Phi_-]$, we obtain 
\begin{align*}
\nr{\widetilde{v}-\widetilde{w}}_{Z[\Phi_-]}
\lesssim
\nr{v-w}_{Z[\Phi_-]}
\braa{
\nr{v}_{Z[\Phi_-]}
\vee
\nr{w}_{Z[\Phi_-]}
}^2,
\end{align*}
where $\widetilde{w}:=w[\Phi_-] +\lambda \mathcal{I}_s[w^3]$. 
The implicit constants in these two estimates are
independent of $\Phi_-$ and $\eta$.
Since $Z$ is complete and convergence in $Z$ implies
convergence in distributions, the continuity of
$\Px^\alpha P^n$ on distributions shows that
$Z[\Phi_-]$ is a Banach space.
Therefore, 
the standard contraction argument yields
a constant $\eta_2\in(0,\eta_1]$ such that 
if $0<\eta\le\eta_2$, then Equation (\ref{Ieq})
has a solution in $Z[\Phi_-]$. 
By the uniqueness in Proposition \ref{prop:S},
this solution coincides with $u$.
The constant $\eta_2$ can be chosen independently of $j$.
Hence we obtain (\ref{PPu}).
\end{proof}

\begin{prop}\label{prop:u-r-d}%\paragraph{prop:u-r-d}
Under the assumptions of Proposition \ref{prop:r-d},
the corresponding solution $u$ satisfies 
\[
\D u(t) \in\Sc\oplus\Sc
\quad\text{and}\quad
U(-t)\D u(t) \in\Sc\oplus\Sc
\]
for any $t\in\R$.
\end{prop}

\begin{proof}
Fix $j\in\{1,2,3\}$ and write $P=P_j$.
Define
\[
H^\infty=\bigcap_{m\in\Z}H^m,
\quad
u_n=P^n u,
\quad
F_n=\lambda P^n(u^3) 
\quad (n\in\Z).
\]
Proposition \ref{prop:r-d} gives 
$u_n\in C(\R;H^{m+1})\cap C^1(\R;H^m)$
for any $m,n\in\Z$.
The Leibniz rule, 
the spatial regularity above,
and the Sobolev embedding $H^2\hookrightarrow L^\infty$
yield
\[
F_n
=
\lambda 
\sum_{m_1+m_2+m_3=n}
\frac{n!}{m_1!m_2!m_3!}
u_{m_1}u_{m_2}u_{m_3}\in C(\R;H^m) 
\quad (n,m\in \Z).
\] 
Furthermore, 
using (\ref{NLKG}) and the spatial regularity above,
we obtain $u\in C^\infty(\R;H^m)$ for any $m\in\Z$
by induction on the number of time derivatives.
The commutation of $P$ with $\Pt^2-\Delta+1$ gives 
\[
\Pt^2u_n=-\omega^2u_n+F_n
\]
in the sense of distributions.
The right-hand side belongs to $C(\R;H^m)$
for any $m\in\Z$, 
which yields
\[
u_n\in C^2(\R;H^m)
\quad (m,n\in\Z).
\]

Fix $t_0\in\R$.
Henceforth, 
all functions and their time derivatives 
are evaluated at $t=t_0$ unless otherwise stated.
We prove by induction on $r\in \Z$ that 
\begin{align*}
x_j^r u_n, x_j^r \Pt u_n \in H^\infty \quad (n\in \Z).
\end{align*} 
The case $r=0$ follows from the regularity established above.
Let $r\in\N$ and assume that
\[
x_j^q u_n,\ x_j^q\Pt u_n\in H^\infty
\quad (0\le q\le r-1,\ n\in\Z).
\]
Since
\[
x_j^r\Pt u_n
=
x_j^{r-1}(u_{n+1}-t_0\Pj u_n),
\]
the induction hypothesis implies that
$x_j^r\Pt u_n\in H^\infty$ for any $n\in\Z$.
Fix $n\in \Z$. 
Differentiating $u_{n+1}=t\Pj u_n+x_j\Pt u_n$
in time before evaluating at $t=t_0$,
and using the equation for $u_n$, 
we obtain 
\begin{align*}
x_j \omega^2 u_n 
=
-\Pt u_{n+1}+t_0\Pj\Pt u_n+\Pj u_n +x_j F_n.
\end{align*} 
Multiplying this identity by $x_j^{r-1}$ and using
\[
[x_j^r,\omega^2]
=
2r x_j^{r-1}\Pj+r(r-1)x_j^{r-2},
\]
we obtain
\begin{align}
\omega^2(x_j^r u_n)
&=
x_j^{r-1}
\braa{
-\Pt u_{n+1}
+t_0\Pj\Pt u_n
+(1-2r)\Pj u_n
}
-r(r-1)x_j^{r-2}u_n
+x_j^rF_n
\nonumber 
\\
&=:
G_{r,n}+x_j^rF_n. \label{id:proo:prop:u-r-d}
\end{align} 
Here, 
the term  $-r(r-1)x_j^{r-2}u_n$ is omitted when $r = 1$.
By the induction hypothesis, we have
$G_{r,n}\in H^\infty$ for any $n\in\Z$.

Fix $N\in\Z$ and choose $\chi\in C_c^\infty(\R^3)$ so that 
$0\le \chi \le 1$, 
$\chi(x)=1$ if $|x|\le 1$, 
and 
$\chi(x)=0$ if $|x|\ge 2$.
For $R\ge 1$, 
set   
\begin{align*}
\chi_R(x)=\chi(x/R) \quad\text{and}\quad 
z_{n,R}=\chi_R x_j^r u_n \quad (0\le n \le N). 
\end{align*} 
Fix $n\in \{ 0,1,\cdots, N \}$. 
Then we have
\begin{align*}
\omega^2 z_{n,R}
&=
\chi_R G_{r,n}
+\lambda 
\sum_{m_1+m_2+m_3=n}
\frac{n!}{m_1!m_2!m_3!}
u_{m_1}u_{m_2}z_{m_3,R}
+E_{n,R},
\end{align*}
where
\begin{align*}
E_{n,R}
:=[\omega^2,\chi_R](x_j^r u_n)
=
-2r(\Pj\chi_R)x_j^{r-1}u_n
-2x_j^r\nabla\chi_R\cdot\nabla u_n
-(\Delta\chi_R)x_j^r u_n.
\end{align*} 
Set $S\ge R$. 
It follows that  
\begin{align}\label{est:proof:prop:u-r-d}
\sum_{n=0}^N
\nr{ z_{n,R}-z_{n,S} }_{H^2} 
&\le 
\sum_{n=0}^N
\nr{ (\chi_R-\chi_S)G_{r,n} }_2
+
\sum_{n=0}^N
\nr{ E_{n,R}- E_{n,S} }_2 
\nonumber
\\
&\quad +
\sum_{n=0}^N
\sum_{m_1+m_2+m_3=n}
\frac{n!}{m_1!m_2!m_3!}
\nr{
u_{m_1}u_{m_2}z_{m_3,R}
-
u_{m_1}u_{m_2}z_{m_3,S}
}_2.
\end{align} 
Since 
$\supp \nabla \chi_R, \supp \Delta \chi_R \subset \{ R\le |x| \le 2R \}$, 
we have 
\begin{align*}
&
\nr{ E_{n,R} }_2 
\\
&\le 
 2rR^{-1} \nr{\Pj\chi   }_\infty \nr{x_j^{r-1}u_n}_2
+4 \nr{\nabla\chi}_\infty \nr{x_j^{r-1} \nabla u_n}_{L^2(\{ |x|\ge R \})}
+2 R^{-1} \nr{\Delta\chi}_\infty \nr{x_j^{r-1} u_n}_2.
\end{align*} 
The induction hypothesis and the Lebesgue dominated convergence theorem 
imply that 
$
\lim_{R\to\infty}
\nr{ E_{n,R} }_2 =0.
$
We see from $G_{r,n} \in H^\infty$ that 
\begin{align}\label{limit:proof:prop:u-r-d}
\lim_{R,S\to\infty}
\sum_{n=0}^N
\nr{ (\chi_R-\chi_S)G_{r,n} }_2
+
\lim_{R,S\to\infty}
\sum_{n=0}^N
\nr{ E_{n,R}- E_{n,S} }_2 
=0.
\end{align} 
Furthermore, 
we obtain 
\begin{align*}
\sum_{n=0}^N
\sum_{m_1+m_2+m_3=n}
&
\frac{n!}{m_1!m_2!m_3!}
\nr{
u_{m_1}u_{m_2}z_{m_3,R}
-
u_{m_1}u_{m_2}z_{m_3,S}
}_2
\\
&\quad 
\le 
3^N
\max_{0\le m \le N}\braa{  \sup_{|x|\ge R} |u_m(x)|^2 }
\sum_{n=0}^N
\nr{ z_{n,R}-z_{n,S} }_2.
\end{align*} 
Since $u_m\in H^2\subset C_0(\R^3)$ for any $0\le m\le N$,
there exists $R_0\ge1$ such that
\[
3^N \max_{0\le m\le N}
\sup_{|x|\ge R_0}|u_m(x)|^2\le\frac12.
\]
Hence, for any $S\ge R\ge R_0$, we obtain
\begin{align*}
\sum_{n=0}^N\nr{z_{n,R}-z_{n,S}}_{H^2}
&\le
2\sum_{n=0}^N\nr{(\chi_R-\chi_S)G_{r,n}}_2
+
2\sum_{n=0}^N\nr{E_{n,R}-E_{n,S}}_2.
\end{align*}
By (\ref{est:proof:prop:u-r-d}) and
(\ref{limit:proof:prop:u-r-d}),
the family $(z_{n,R})_{R\ge1}$ is Cauchy in $H^2$.
On the other hand,
$(z_{n,R})_{R\ge1}$ converges to $x_j^r u_n$
in the sense of distributions as $R\to\infty$.
Consequently,
\begin{align}\label{xj-H2}
x_j^r u_n\in H^2 \quad (n=0,1,\cdots,N).
\end{align} 
Note that the identity 
\begin{align*} 
\omega^2 (x_j^r u_n) 
=
G_{r,n}
+ \lambda 
\sum_{m_1+m_2+m_3=n}
\frac{n!}{m_1!m_2!m_3!}
u_{m_1}u_{m_2} (x_j^r u_{m_3})
\end{align*}
implies $x_j^r u_n \in H^\infty$. 
Indeed, by the Leibniz rule and the Sobolev embedding,
if $x_j^r u_n\in H^k$ ($0\le n\le N$) for some integer $k\ge 2$,
then the right-hand side belongs to $H^k$
for any $0\le n\le N$.
Hence $x_j^r u_n\in H^{k+2}$ for any $0\le n\le N$.
Starting from (\ref{xj-H2}) and iterating this argument,
we obtain $x_j^r u_n\in H^\infty$
for any $0\le n\le N$.
Since $N\in\Z$ is arbitrary,
this holds for any $n\in\Z$.
Together with the regularity of $x_j^r\Pt u_n$
established above, this completes the induction on $r$.

Taking $n=0$ and recalling that $j$ is arbitrary, we obtain
\[
x_j^r u(t_0),\ x_j^r \Pt u(t_0)\in H^\infty
\quad (r\in\Z,\ j=1,2,3).
\]
By the Leibniz rule and the Sobolev embedding,
these properties imply
\[
\D u(t_0)\in\Sc\oplus\Sc.
\]
Since $t_0\in\R$ is arbitrary,
this holds for any $t\in\R$.
Finally, using the fact that $U(-t)$ preserves
$\Sc\oplus\Sc$ for each fixed $t\in\R$, we obtain
\[
U(-t)\D u(t)\in\Sc\oplus\Sc
\quad (t\in\R).
\]
\end{proof}

\begin{prop}\label{prop:u-decay}%\paragraph{prop:u-decay}
Under the assumptions of Proposition \ref{prop:r-d},
the corresponding solution $u$ satisfies the following properties: 
\begin{enumerate}[(1)]
  \item 
$
\nr{u(t)}_\infty \lesssim t^{-3/2} \quad (t\ge 1).
$ 
  \item 
$
\nr{U(-t)\D u(t)-\s \Phi_-}_{\mh^1}=O(t^{-2})
$  
as $t\to\infty$.
\item
$
\nr{\D P_j^3\braa{u-w[\s\Phi_-]}(t)}_{\mh^{1/2}}
=o(1)
$
as $t\to\infty$ $(j=1,2,3)$.
\end{enumerate}
\end{prop}

\begin{proof}
\begin{enumerate}[(1)]
  \item 
By Proposition \ref{prop:r-d}, 
we have 
$M:=\sup_{t\ge 0}\nr{u(t)}_{H^3}<\infty$ 
and 
$u\in L^4 W^1_4$. 
The embedding $W^1_4 \hookrightarrow L^\infty$ implies that 
\[
\delta_T:=\nr{u}_{L^4([T,\infty);L^\infty)}\to 0 
\quad\text{as $T\to\infty$.}
\]
For $T\ge 1$ and $R>T$, 
we set 
\[
D_{T,R}=\sup_{T\le t \le R} t^{3/2}\nr{u(t)}_\infty,
\] 
which is finite by the continuity of $u$ in $H^2$.
Note that the solution $u$ satisfies 
\begin{align}\label{eq:proof:prop:u-decay} 
u(t)=\braa{ U(t-T)\D u(T) }_1 + \lambda 
\int_T^t \omega^{-1}\sin((t-\tp)\omega)u(\tp)^3 d\tp 
\quad (t\ge T).
\end{align}
Using Proposition \ref{prop:free-decay} and
the embedding $H^2\hookrightarrow L^\infty$,
we obtain
\begin{align*}
&
\nr{\braa{U(t-T)\D u(T)}_1}_\infty
\\
&\lesssim
\brab{
|t-T|^{-3/2}
\braa{
\nr{u(T)}_{W^3_1}
+\nr{\Pt u(T)}_{W^2_1}
}
}
\wedge
\braa{
\nr{u(T)}_{H^2}
+\nr{\Pt u(T)}_{H^1}
}
\quad (t>T).
\end{align*}
By Proposition \ref{prop:u-r-d}, 
we have 
\begin{align*}
A_T:=
\sup_{t\ge T} t^{3/2} \nr{ \braa{ U(t-T)\D u(T) }_1 }_\infty <\infty.
\end{align*} 
Set $L\ge 1$ and $T\ge 2L$, 
and define for $t\ge T$, 
\[
I_1=[T,t] \cap (-\infty,t/2],
\quad
I_2=[T,t] \cap [t/2,t-L],
\quad
I_3=[T,t] \cap [t-L,t].
\]
By Proposition \ref{prop:free-decay}, the Leibniz rule,
H\"older's inequality and the Gagliardo--Nirenberg inequality,
we obtain 
\begin{align*}
&
t^{3/2}\int_{I_1} \nr{\omega^{-1}\sin((t-\tp)\omega)u(\tp)^3}_\infty d\tp 
\lesssim
t^{3/2}\int_{I_1} |t-\tp|^{-3/2} \nr{u(\tp)^3}_{W^3_1} d\tp 
\\
&\lesssim
t^{3/2}
\int_{I_1} t^{-3/2} \nr{u(\tp)}_{H^3}^2 \nr{u(\tp)}_\infty d\tp 
\\
&\lesssim
\int_{I_1} (\tp)^{-3/2} (\tp)^{3/2}\nr{u(\tp)}_\infty 
\nr{u(\tp)}_{H^3}^2 d\tp 
\\
&\lesssim
D_{T,R}
\int_T^\infty (\tp)^{-3/2}  
\nr{u(\tp)}_{H^3}^2 d\tp 
\lesssim
D_{T,R}M^2 T^{-1/2} \quad (T\le t \le R).
\end{align*} 
Similarly, we have
\begin{align*}
&
t^{3/2}\int_{I_2}
\nr{\omega^{-1}\sin((t-\tp)\omega)u(\tp)^3}_\infty d\tp
\lesssim
M^2D_{T,R}
\int_{I_2}
(t-\tp)^{-3/2}
\braa{\frac{t}{\tp}}^{3/2} d\tp
\\
&\lesssim
M^2D_{T,R}
\int_L^\infty \tau^{-3/2}\,d\tau
\lesssim
M^2L^{-1/2}D_{T,R}
\quad (T\le t\le R).
\end{align*}
We also see from the embedding $H^2 \hookrightarrow L^\infty$ 
and H\"older's inequality in time that 
\begin{align*}
&
t^{3/2}\int_{I_3}
\nr{\omega^{-1}\sin((t-\tp)\omega)u(\tp)^3}_\infty d\tp
\lesssim
t^{3/2}\int_{I_3}
\nr{u(\tp)^3}_{H^1} d\tp
\\
&\lesssim
t^{3/2}\int_{I_3}
\nr{u(\tp)}_{H^1} \nr{u(\tp)}_\infty^2 d\tp
\\
&\lesssim
M D_{T,R}
\int_{I_3} \nr{u(\tp)}_\infty d\tp
\\
&\lesssim
M D_{T,R} \delta_T |I_3|^{3/4}
\lesssim
M L^{3/4}D_{T,R} \delta_T 
\quad (T\le t\le R).
\end{align*}
It follows from (\ref{eq:proof:prop:u-decay}) that 
\begin{align*}
D_{T,R} \le A_T+C(M^2T^{-1/2}+M^2L^{-1/2}+ML^{3/4}\delta_T)D_{T,R}.
\end{align*} 
Here, 
$C$ is a positive constant independent of $M,T,L$ and $R$. 
Let $L\ge 1$ satisfy $CM^2L^{-1/2}\le 1/4$. 
We next choose $T\ge 2L$ so that 
\[
C\braa{M^2T^{-1/2}+ML^{3/4}\delta_T}\le 1/4.
\]
Then we obtain $D_{T,R}\le 2A_T$ for any $R>T$.
As $R\to\infty$, 
we have 
\[
\sup_{t\ge T}t^{3/2}\nr{u(t)}_\infty\le2A_T<\infty.
\]
The continuity of $u$ in $H^2$
on $[1,T]$ implies the desired estimate.
  \item
By (1) and the uniform boundedness
of $\nr{u(t)}_2$, we have
\[
\nr{u(t)^3}_2
\le \nr{u(t)}_\infty^2\nr{u(t)}_2
\lesssim t^{-3}
\quad (t\ge 1).
\]
By (\ref{conv}) and the Duhamel formula, 
we obtain 
\[
U(-t)\D u(t)-\s\Phi_-
=
-\lambda\int_t^\infty
U(-\tp)\binom{0}{u(\tp)^3}\,d\tp
\]
in $\mh^{1/2}$.
For $t\ge 1$, the integral on the right also converges
absolutely in $\mh^1$, 
since
$\| u(\tp)^3 \|_2 \lesssim (\tp)^{-3}$.
Hence, by the unitarity of $U(t)$ on $\mh^1$,
\[
\nr{U(-t)\D u(t)-\s\Phi_-}_{\mh^1}
\le
\int_t^\infty\nr{u(\tp)^3}_2\,d\tp
\lesssim t^{-2}
\qquad (t\ge 1).
\]
  \item 
Write $P=P_j$ for a fixed $j\in\{1,2,3\}$, and set
\[
v=P^3(u-w[\s\Phi_-]),
\qquad F=\lambda u^3.
\]
By the future Duhamel representation and
the future version of Proposition \ref{prop:PI}, we have
\[
\D v(t)
=
-\int_t^\infty
U(t-\tp)\binom{0}{P^3F(\tp)}\,d\tp.
\]
Hence Proposition \ref{prop:St} implies 
\[
\nr{\D v(t)}_{\mh^{1/2}}
\lesssim
\nr{P^3F}_{L^{4/3}((t,\infty);L^{4/3})}.
\]
By the Leibniz rule and H\"older's inequality,
the right-hand side is bounded by
\[
\sum_{m_1+m_2+m_3=3}
\frac{3!}{m_1!m_2!m_3!}
\prod_{\ell=1}^3
\nr{P^{m_\ell}u}_{L^4((t,\infty);L^4)}.
\]
Since $P^m u\in Z\subset L^4L^4$ for any
$m\in\{0,1,2,3\}$ by Proposition \ref{prop:r-d},
this expression tends to zero as $t\to\infty$.
\end{enumerate}
\end{proof}

\begin{prop}\label{prop:weighted-growth}%\paragraph{prop:weighted-growth}
Under the assumptions of Proposition \ref{prop:r-d},
the corresponding solution $u$ satisfies
\begin{align}\label{est:prop:weighted-growth}
\nr{x_j^r P_j^n u(t)}_{H^1}
+
\nr{x_j^r \Pt P_j^n u(t)}_2
\le C t^r \quad (t\ge T)
\end{align} 
for any $r,n\in\Z$ with $r+n\le3$,
and any $j\in\{1,2,3\}$. 
Here, 
$C>0$ and $T\ge 1$ are constants 
independent of $r,n$ and $t$. 
\end{prop}

\begin{proof}
Fix $j\in\{1,2,3\}$,
and write $P=P_j$, $u_n=P^n u$ and
$F_n=\lambda P^n(u^3)$.
By Proposition \ref{prop:u-decay}(1),
there exists $T_0\ge 1$ such that
$\nr{u(t)}_\infty\le 1/2$ for any $t\ge T_0$. 
Therefore, 
if $t\ge T_0$, 
then for any $f\in H^1$
\begin{align*}
\nr{f}_{H^1}^2 
&= \Re \bra{(\omega^2-\lambda u(t)^2) f,f}_{H^{-1},H^1} 
+ \lambda \, \Re \int_{\R^3} u(t)^2 |f|^2 dx 
\\
&\le 
\frac{1}{4}\nr{f}_{H^1}^2
+\nr{f}_{H^1} \nr{(\omega^2-\lambda u(t)^2)f}_{H^{-1}}, 
\end{align*} 
which implies that 
\begin{align}\label{est:proof:prop:weighted-growth}
\nr{f}_{H^1}
\le 
2\nr{(\omega^2-\lambda u(t)^2)f}_{H^{-1}}.
\end{align}

We prove (\ref{est:prop:weighted-growth}) by induction on $r$.
The case $r=0$ follows from Proposition \ref{prop:r-d}. 
Let $k\in\{1,2,3\}$ and assume that
(\ref{est:prop:weighted-growth}) holds
for any $r,n\in\Z$ with
$0\le r \le k-1$ and $r+n\le 3$. 
Let $T\ge 1$ be a common time beyond which
these estimates hold.
We fix $n\in \{ 0,\cdots,3-k\}$.
Using the identities  
$x_j^k \Pt u_n=x_j^{k-1}(u_{n+1}-t\Pj u_n)$ 
and 
$x_j^{k-1} \Pj u_n=\Pj(x_j^{k-1}u_n)-(k-1)x_j^{k-2}u_n$, 
we obtain
$\nr{x_j^k \Pt u_n}_2\lesssim t^k$ for any $t\ge T$. 
Recall the identity (\ref{id:proo:prop:u-r-d})
in the proof of Proposition \ref{prop:u-r-d},
with $r=k$ and $t_0=t$.
The induction hypothesis implies that
$\nr{G_{k,n}}_{H^{-1}}\lesssim t^k$ for any $t\ge T$.
We see from (\ref{est:proof:prop:weighted-growth}) that 
\begin{align*}
\nr{ x_j^k u }_{H^1} 
\le 2 \nr{ (\omega^2-\lambda u^2)(x_j^k u) }_{H^{-1}}
=   2 \nr{ G_{k,0} }_{H^{-1}} 
\lesssim t^k \quad (t\ge T \vee T_0).
\end{align*} 
If we set $a_1=3\lambda u u_1$ and $a_2=3\lambda uu_2+6\lambda u_1^2$, 
then $F_n=a_n u$ when $n=1,2$.
Since $a_1, a_2 \in C_b H^2 \subset C_b L^\infty$, 
we have for any $n=1,2$ with $k+n\le 3$, 
\begin{align*}
\nr{x_j^k u_n}_{H^1}&=\nr{\omega^2 (x_j^k u_n)}_{H^{-1}} 
\le \nr{G_{k,n}}_{H^{-1}}+\nr{a_n (x_j^k u)}_{H^{-1}} 
\\
&\le \nr{G_{k,n}}_{H^{-1}}+\nr{a_n}_\infty \nr{ x_j^k u }_2 
\lesssim t^k \quad (t\ge T \vee T_0). 
\end{align*}  
Together with the estimate for $x_j^k\Pt u_n$,
this completes the induction and proves the proposition. 
\end{proof}

\section{Construction of test functions}
The goal of this section is to prove the following lemma:
\begin{lem}\label{lem:main}%\paragraph{lem:main}
Under the assumptions of Proposition \ref{prop:r-d},
suppose further that $\Phi_-$ is real-valued and non-zero. 
Let $u$ be the solution to (\ref{Ieq}). 
Then, for any $j=1,2,3$, 
there exist a constant $C_\psi>0$ and a family of functions
$(\psi_t)_{t\ge1}\subset C_c^\infty(\R^3)$
satisfying (\ref{property:psi:1})--(\ref{property:psi:3}).
\end{lem}

\begin{proof}
We divide the proof into four steps.
\\
(Step 1) %\paragraph{step 1}
Write $\s \Phi_-={}^t (f_+,g_+)$ and $\phi=(f_+ -i\omega^{-1}g_+)/2$. 
Then $u_+(t)=2\Re(e^{it\omega}\phi)$.
The injectivity of $\s$ and the identity $\s0=0$
imply that $\s\Phi_-\ne0$.
Since $\s\Phi_-$ is real-valued, 
it follows that $\phi\ne0$.
Proposition \ref{prop:r-d} also gives $\phi\in\Sc$.
Define $v(t,x)=(e^{it\omega}\phi)(x)$.
By the stationary phase method
(see, e.g., \cite[Section 7.7]{Hormander2003}), 
for any compact subset
$K\subset \{y\in\R^3\,;\,|y|<1\}$,
there exists a constant $C_K>0$ such that
\begin{align}\label{MSP} 
\sup_{y\in K}
\left|
v(t,ty)-\frac{1}{2}t^{-3/2}a(y)e^{it\rho(y)}
\right|
\le C_Kt^{-5/2}
\quad (t\ge1).
\end{align}
Here, 
\[
\rho(y)=\sqrt{1-|y|^2},
\qquad
a(y)=
2
e^{3\pi i/4}\rho(y)^{-5/2}
\widehat{\phi}\left(-\frac{y}{\rho(y)}\right).
\]
Fix $j\in \{ 1,2,3 \}$. 
Since $\phi\neq0$, the function $a$ is not identically zero.
Hence we can choose a closed ball $K_1$ of positive radius 
such that
\[
K_1 \subset \{y\in\R^3\,;\,|y|<1, \ y_j\neq 0 \},
\qquad
a(y)\neq 0 \quad (y\in K_1),
\]
and such that a smooth branch of $\arg a$
exists on a neighborhood of $K_1$. 
Set $\theta_t(y)=t\rho(y)+\arg a(y)$
on a neighborhood of $K_1$. 
Then we have 
\[
\sup_{y\in K_1}
\left|
u_+(t,ty)-t^{-3/2}|a(y)|\cos \theta_t(y)
\right|
\le 2C_{K_1}t^{-5/2}
\quad (t\ge1).
\]
Applying the same stationary phase argument to $\Pj v$,
and enlarging $C_{K_1}$ if necessary, we also obtain
\begin{align*}
\sup_{y\in K_1}
\left|
(\Pj u_+)(t,ty)
-t^{-3/2}\frac{y_j}{\rho(y)}
|a(y)|\sin\theta_t(y)
\right|
\le 2C_{K_1}t^{-5/2}
\quad (t\ge1).
\end{align*}
Since 
\[
\sup_{t\ge 1}
\sup_{y\in K_1}
\braa{
\abs{
|a(y)|\cos \theta_t(y)
}
+
\abs{
\frac{y_j}{\rho(y)}
|a(y)|\sin\theta_t(y)
}}<\infty,
\]
we have 
\begin{align}\label{est:proof:lem:main}
\sup_{y\in K_1}
\left|
(\Pj F_+)(t,ty)
-3\lambda t^{-9/2}\frac{y_j}{\rho(y)}
|a(y)|^3\cos^2\theta_t(y)\sin\theta_t(y)
\right|
\le C_1t^{-11/2}
\quad (t\ge1).
\end{align} 
Here, 
$C_1$ is a positive constant independent of $t$.
\\ \ \\
(Step 2) %\paragraph{step 2}
Choose a real-valued function
$\chi\in C_c^\infty(\R^3)$ such that 
$\supp \chi \subset K_1^\circ$, 
$0\le\chi\le1$ and $\chi\not\equiv 0$.
Define
$\chi_t(x)=\chi(x/t)$
for $t\ge 1$ and $x\in\R^3$.
For $t\ge1$, we also define $b_t$ by
\[
b_t(ty)
=
3\lambda t^{-3/2}\chi(y)\frac{y_j^3}{\rho(y)}
|a(y)|^3\cos^2\theta_t(y)\sin\theta_t(y)
\quad (y\in K_1),
\]
and set $b_t(ty)=0$ for $y\notin K_1$.
By construction, $b_t\in C_c^\infty(\R^3)$.
Estimate (\ref{est:proof:lem:main}) gives 
\begin{align}\label{b+e}
\chi_t (x) t x_j^2 \Pj F_+(t,x)
=
b_t(x)+e_t(x) \quad (x\in \R^3),
\end{align} 
where $e_t\in L^2$ satisfies
$\nr{e_t}_2=O(t^{-1})$ as $t\to\infty$.
Indeed, by the change of variables $x=ty$, we have
\[
\nr{e_t}_2
\le
C_1t^{-1}
\left(
\int_{K_1}\chi(y)^2y_j^4\,dy
\right)^{1/2}.
\]
We define $\psi_t=\chi_t b_t$ for $t\ge1$.
Since
\[
\Px\bigl(\theta_t(x/t)\bigr)
=
(\Px\rho)(x/t)
+t^{-1}(\Px\arg a)(x/t)
\]
on a neighborhood of $\supp\psi_t$,
the chain rule yields
\[
\nr{\Px^\alpha\psi_t}_\infty
\lesssim t^{-3/2}
\quad (t\ge1,\ |\alpha|\le2).
\]
Together with
$\supp\psi_t\subset t\,\supp\chi$,
this implies (\ref{property:psi:1}).
Furthermore, by the definition of $b_t$ and
the change of variables $x=ty$, we have
\begin{align} 
\nr{b_t}_2^2
&=9
\int_{K_1} \chi(y)^2 \frac{y_j^6}{\rho(y)^2}
|a(y)|^6 \cos^4\theta_t(y)\sin^2\theta_t(y)
dy
\label{id:bt}
\\
&\le
9\int_{K_1}
\chi(y)^2\frac{y_j^6}{\rho(y)^2}|a(y)|^6\,dy
<\infty
\quad (t\ge1).
\label{id:bt:2}
\end{align}
Hence, by (\ref{b+e}), we obtain 
\begin{align}\label{bra}
\bra{tx_j^2 \Pj F_+,\psi_t}=\nr{b_t}_2^2+\bra{e_t,b_t }
=\nr{b_t}_2^2+O(t^{-1})
\quad\text{as $t\to\infty$.}
\end{align} 
\ \\
(Step 3) %\paragraph{step 3}
If we prove 
\begin{align}\label{final-task}
\nr{b_t}_2^2
=
\frac{9}{16}
\int_{K_1}
\chi(y)^2\frac{y_j^6|a(y)|^6}{\rho(y)^2}\,dy 
+O(t^{-1}) \quad \text{as $t\to\infty$,}
\end{align}
then (\ref{bra}) implies (\ref{property:psi:3}) with
\[
C_\psi:=
\frac{9}{16}
\int_{K_1}
\chi(y)^2\frac{y_j^6|a(y)|^6}{\rho(y)^2}\,dy>0.
\]
Define 
$
A(y)=\chi(y)^2 y_j^6|a(y)|^6/\rho(y)^2
$
($y\in K_1$). 
For any non-zero integer $\ell$, integration by parts gives
\begin{align*}
&
\int_{K_1}
e^{i\ell t\rho(y)}
A(y)e^{i\ell\arg a(y)}\,dy
\\
&=
-\frac{1}{i\ell t}
\int_{K_1}
e^{i\ell t\rho(y)}
\Pj\left(
\frac{A(y)e^{i\ell\arg a(y)}}{\Pj\rho(y)}
\right)\,dy
=O(t^{-1})
\quad\text{as $t\to\infty$}.
\end{align*}
By (\ref{id:bt}) and the identity
\[
\cos^4\theta\sin^2\theta
=
\frac{1}{16}
+\frac{1}{32}\cos(2\theta)
-\frac{1}{16}\cos(4\theta)
-\frac{1}{32}\cos(6\theta),
\]
we obtain (\ref{final-task}).
\\ \ \\
(Step 4) %\paragraph{step 4}
Let us show 
(\ref{property:psi:2}). 
Write $P=P_j$.
By the definition of $d(t)$,
the Cauchy--Schwarz inequality, 
integration by parts and (\ref{id:bt:2}),  
we have 
\begin{align}\label{est:d-psi}
\abs{ 
\bra{d(t),\psi_t}_{L^2} 
}
&\le 
\braa{
\nr{P^3(u-u_+)(t)}_2  
+
3\nr{\chi_t tx_jF(t)}_2 
+ 
 \nr{\chi_t x_j^2 PF(t)}_2
}
\nr{b_t}_2
\nonumber
\\
&\quad + 
\abs{
\bra{2tx_j^2\Pj (F(t)-F_+(t)) ,\psi_t}_{L^2} 
}
\nonumber
\\
&\lesssim 
\nr{P^3(u-u_+)(t)}_2  
+
\nr{\chi_t tx_jF(t)}_2 
+ 
\nr{\chi_t x_j^2 PF(t)}_2
\nonumber
\\
&\quad +
\abs{
\bra{ F(t)-F_+(t), t\Pj(x_j^2\psi_t)}_{L^2} 
}
\end{align} 
for any $t\ge 1$. 
Propositions \ref{prop:u-decay}(1) and \ref{prop:r-d} imply that 
\begin{align*}
\nr{\chi_t tx_jF(t)}_2 
\le 
t^2 \nr{u(t)^3}_{L^2(\{ |x|\le  t \})} 
\le 
t^2 \nr{u(t)}_\infty^2 \nr{u(t)}_2 
\lesssim
t^{-1} 
\end{align*} 
and 
\begin{align*}
\nr{\chi_t x_j^2 PF(t)}_2 
\le 
3t^2 \nr{u(t)^2 Pu(t)}_{L^2(\{ |x|\le  t \})} 
\le 
3t^2 \nr{u(t)}_\infty^2 
\nr{Pu(t)}_2  
\lesssim
t^{-1} 
\end{align*} 
for any $t\ge 1$. 
The uniform $H^2$ bound established in Step 2 also gives 
\begin{align*}
\nr{t\Pj(x_j^2\psi_t)}_2
&\le
2\nr{tx_j\psi_t}_2
+\nr{tx_j^2\Pj\psi_t}_2
\le
2t^2\nr{\psi_t}_2
+t^3\nr{\Pj\psi_t}_2
\lesssim t^3
\quad (t\ge 1).
\end{align*}
Propositions \ref{prop:u-decay}
and \ref{prop:free-decay} yield 
\begin{align*}
\nr{F(t)-F_+(t)}_2 
&\lesssim
\nr{u(t)-u_+(t)}_2 
\braa{ \nr{u(t)}_\infty \vee \nr{u_+(t)}_\infty }^2
\\
&\lesssim
t^{-3} \nr{u(t)-u_+(t)}_2 
=
O(t^{-5}) 
\quad \text{as $t\to\infty$.}
\end{align*} 
It follows from (\ref{est:d-psi}) and Proposition \ref{prop:u-decay}(3) 
that (\ref{property:psi:2}) holds true. 
\end{proof}

\section{Proof of the main theorem}
In this section, we prove Theorem \ref{thm:main}.
Let $\eta_2$ be the positive number appearing 
in Proposition \ref{prop:r-d}. 
Assume that $\Phi_- \in B_{\eta_2}\mh^{1/2} \cap (\Sc \oplus \Sc)$ 
is non-zero and real-valued.  
Let $u$ be the corresponding solution to (\ref{Ieq}). 
By Propositions \ref{prop:r-d} and \ref{prop:u-r-d},
we have $\s\Phi_-\in\Sc\oplus\Sc$ and
$U(-t)\D u(t)\in\Sc\oplus\Sc$ for every $t\in\R$.
Set $u_+=w[\s \Phi_-]$, $F=\lambda u^3$ and $F_+=\lambda u_+^3$. 
Then we see from 
Proposition \ref{prop:AUUB} that for any $j=1,2,3$, 
\begin{align}\label{id:DP3}
&
\D P_j^3(u-u_+)
\nonumber
\\
&=A_j(t)^3\D(u-u_+)
+
 A_j(t)^2\binom{0}{x_jF}
+A_j(t)\binom{0}{x_jP_jF}
+\binom{0}{x_jP_j^2F}
\nonumber
\\
&=A_j(t)^3\D(u-u_+)
+
\binom{
3tx_jF+2tx_j^2\Pj F+x_j^2P_jF
}{
(\Pj-x_j\omega^2)(x_j^2F)
+t^2\Pj^2(x_jF)
+t\Pj(x_jP_jF)
+x_jP_j^2F
}. 
\end{align}
We first show the third inequality in (\ref{core}). 
Combining Propositions \ref{prop:B-bdd} and \ref{prop:AUUB}
with the unitarity of $U(t)$ on $\mh^s$ ($s\in\R$),
we obtain 
\begin{align*}
\nr{U(-t)\D u(t)}_{\mh^{1,3}} 
&\lesssim
\nr{U(-t)\D u(t)}_{\mh^1}
+
\sum_{j=1}^3 \nr{B_j^3 U(-t)\D u(t) }_{\mh^{-2}}
\\
&=\nr{U(-t)\D u(t)}_{\mh^1}
+
\sum_{j=1}^3 \nr{U(-t)A_j(t)^3 \D u(t) }_{\mh^{-2}}
\\
&=\nr{\D u(t)}_{\mh^1}
+
\sum_{j=1}^3 \nr{A_j(t)^3 \D u(t) }_{\mh^{-2}}
\quad (t\in \R).
\end{align*} 
By Proposition \ref{prop:r-d}, we have
\[
\sup_{t\in\R}
\left(
\nr{\D u(t)}_{\mh^1}
+
\sum_{j=1}^3
\nr{\D P_j^3u(t)}_{\mh^{-2}}
\right)<\infty.
\]
We see from (\ref{id:DP3}) and the identity 
$A_j(t)^3\D u_+(t)=\D P_j^3u_+(t)$ that 
\begin{align*}
\nr{U(-t)\D u(t)}_{\mh^{1,3}}
\lesssim
1+\sum_{j=1}^3 R_j(t) \quad (t\in \R).
\end{align*} 
Here, 
for $j=1,2,3$, 
we set
\begin{align}
R_j(t)
&=
\nr{
3tx_jF+2tx_j^2\Pj F+x_j^2P_jF
}_{H^{-2}}
\nonumber 
\\
&\quad+
\nr{
(\Pj-x_j\omega^2)(x_j^2F)
+t^2\Pj^2(x_jF)
+t\Pj(x_jP_jF)
+x_jP_j^2F
}_{H^{-3}}  \quad (t\in \R).
\label{def:Rj}
\end{align}
Using Propositions \ref{prop:u-decay}(1) and \ref{prop:weighted-growth}, 
the embedding $H^2 \hookrightarrow L^\infty$, 
and the uniform $H^2$ bound for $P_j u$ 
provide by Proposition \ref{prop:r-d}, 
we obtain 
\begin{align*}
\nr{x_j P_j^2 F(t)}_2 
&\lesssim
\nr{x_j u(t)^2 P_j^2 u(t) }_2 + \nr{x_j u(t) \braa{ P_j u(t)}^2 }_2
\\
&\lesssim 
\nr{x_j P_j^2 u(t) }_2 \nr{u(t)}_\infty^2  + 
\nr{x_j P_j u(t) }_2 \nr{P_j u(t) }_\infty \nr{u(t)}_\infty 
\\
&\lesssim
t^{1-3}+t^{1-3/2}
\lesssim 
t^{-1/2} 
\end{align*}  
for all sufficiently large $t>0$.
Similarly, 
we have 
\begin{align*}
\nr{x_j^r F(t)}_2 =O(t^{r-3})
\quad (r=0,1,2,3)  
\end{align*} 
and 
\begin{align*}
\nr{x_j^r P_j F(t)}_2 =O(t^{r-3}) 
\quad (r=0,1,2)   
\end{align*} 
as $t \to\infty$.
It follows from (\ref{def:Rj}) that 
\begin{align*}
R_j(t)
&\lesssim
t\|x_jF\|_2+t\|x_j^2F\|_2+\|x_j^2P_jF\|_2
+\|x_j^3F\|_2+\|x_j^2F\|_2+\|x_jF\|_2 
\\
&\quad +
 t^2\|x_jF\|_2+t\|x_jP_jF\|_2
+\|x_jP_j^2F\|_2
=O(1) \quad\text{as $t\to\infty$}. 
\end{align*} 
Since $\s\Phi_-\in\Sc\oplus\Sc\subset\mh^{1,3}$,
the triangle inequality yields the third inequality
in (\ref{core}).
\\

We next show the first inequality in (\ref{core}).
Fix $j\in \{ 1,2,3 \}$ 
and write $P=P_j$, $A(t)=A_j(t)$ and $B=B_j$. 
We see from 
Propositions \ref{prop:B-bdd} and \ref{prop:AUUB} 
and 
the unitarity of $U(t)$ on $\mh^{-2}$ that for any $t\in \R$, 
\begin{align}\label{est:1:h(t)}
\nr{U(-t)\D u(t)-\s \Phi_-}_{\mh^{1,3}} 
&\gtrsim 
\nr{B^3 \braa{ U(-t)\D u(t)-\s \Phi_-}}_{\mh^{-2}} 
\nonumber
\\
&=
\nr{U(t) B^3 U(-t) \braa{ \D u(t)-U(t) \s \Phi_-}}_{\mh^{-2}} 
\nonumber
\\
&=
\nr{A(t)^3 \D \braa{ u(t)-u_+(t) }}_{\mh^{-2}}  
\nonumber
\\
&\ge 
\nr{\braa{ A(t)^3 \D \braa{ u(t)-u_+(t) } }_1 }_{H^{-2}}.
\end{align} 
By (\ref{id:DP3}), 
we have 
\begin{align*}
2tx_j^2 \Pj F_+
=
d(t)-
\braa{ A(t)^3 \D \braa{ u(t)-u_+(t) } }_1 .
\end{align*} 
Here, 
$d(t)$ is the function defined in (\ref{def:d}).
Using 
(\ref{est:1:h(t)}) and 
the family of functions $(\psi_t)_{t\ge 1}$ 
appearing in Lemma \ref{lem:main}, 
we have 
\begin{align*}
2C_\psi 
&=
\abs{ \bra{ 2tx_j^2 \Pj F_+,\psi_t }_{L^2} }+O(t^{-1})
\\
&\le 
\abs{ \bra{ d(t),\psi_t }_{L^2} }
+
\nr{\braa{ A(t)^3 \D \braa{ u(t)-u_+(t) } }_1 }_{H^{-2}} 
\nr{\psi_t}_{H^2}+O(t^{-1})
\\
&\lesssim
\nr{\braa{ A(t)^3 \D \braa{ u(t)-u_+(t) } }_1 }_{H^{-2}}
+o(1)
\\
&\lesssim
\nr{U(-t)\D u(t)-\s \Phi_-}_{\mh^{1,3}}
+
o(1) \quad \text{as $t\to \infty$}.
\end{align*} 
Since $C_\psi>0$, 
we obtain 
the first inequality in (\ref{core}). 
\\

Recall that
\[
h(t)=U(-t)\D u(t)-\s\Phi_-.
\]
By Proposition \ref{prop:u-decay}(2), we have
\[
\lim_{t\to+\infty}\nr{h(t)}_{\mh^1}=0.
\]
If $0\le\sigma<3$, then weighted interpolation gives
\[
\nr{h(t)}_{\mh^{1,\sigma}}
\le
\nr{h(t)}_{\mh^1}^{1-\sigma/3}
\nr{h(t)}_{\mh^{1,3}}^{\sigma/3}.
\]
By the third inequality in (\ref{core}),
$\nr{h(t)}_{\mh^{1,3}}$ is bounded
for all sufficiently large $t$.
Since $\nr{h(t)}_{\mh^1}\to0$ as $t\to+\infty$
and $1-\sigma/3>0$, we obtain
\[
\lim_{t\to+\infty}\nr{h(t)}_{\mh^{1,\sigma}}=0,
\]
which proves (\ref{conv:1}).

We next assume that $\sigma>3$.
By weighted interpolation, 
we have 
\[
\nr{h(t)}_{\mh^{1,3}}
\le
\nr{h(t)}_{\mh^1}^{1-3/\sigma}
\nr{h(t)}_{\mh^{1,\sigma}}^{3/\sigma}.
\]
We see from 
the first inequality in (\ref{core}) 
that 
there exists a constant $c>0$ such that
\[
c\le
\nr{h(t)}_{\mh^1}^{1-3/\sigma}
\nr{h(t)}_{\mh^{1,\sigma}}^{3/\sigma}
\]
for all sufficiently large $t$.
Since $1-3/\sigma>0$ and
$\nr{h(t)}_{\mh^1}\to0$ as $t\to+\infty$,
we obtain
\[
\lim_{t\to+\infty}\nr{h(t)}_{\mh^{1,\sigma}}=\infty.
\]
This proves (\ref{conv:2})
and completes the proof of Theorem \ref{thm:main}.

\section*{Acknowledgments}
The author used ChatGPT (OpenAI) for assistance with English editing, manuscript organization, 
the development and refinement of proof arguments, 
and the review of mathematical arguments and references. 
All mathematical arguments and references were checked by the author, 
who reviewed and revised the manuscript
 and takes full responsibility for its content.

\bib
\end{document}